
\magnification=1100
\overfullrule0pt

\input amssym.def
\input prepictex
\input pictex
\input postpictex



\def\CC{{\Bbb C}}
\def\FF{{\Bbb F}}

\def\ZZ{{\Bbb Z}}

\def\fg{{\frak{g}}}
\def\fgl{{\frak{gl}}}

\def\deg{\hbox{deg}}
\def\End{\hbox{End}}

\def\id{\hbox{id}}



\font\smallcaps=cmcsc10
\font\titlefont=cmr10 scaled \magstep1

\font\sectionfont=cmbx10
\font\tinyrm=cmr10 at 8pt


\newcount\sectno
\newcount\subsectno
\newcount\resultno

\def\section #1. #2\par{
\sectno=#1
\resultno=0
\bigskip\noindent{\sectionfont #1.  #2}~\medbreak}

\def\subsection #1\par{\bigskip\noindent{\it  #1} \medbreak}


\def\prop{ \global\advance\resultno by 1
\bigskip\noindent{\bf Proposition \the\sectno.\the\resultno. }\sl}
\def\lemma{ \global\advance\resultno by 1
\bigskip\noindent{\bf Lemma \the\sectno.\the\resultno. }
\sl}
\def\remark{ \global\advance\resultno by 1
\bigskip\noindent{\bf Remark \the\sectno.\the\resultno. }}
\def\example{ \global\advance\resultno by 1
\bigskip\noindent{\bf Example \the\sectno.\the\resultno. }\sl}
\def\cor{ \global\advance\resultno by 1
\bigskip\noindent{\bf Corollary \the\sectno.\the\resultno. }\sl}
\def\thm{ \global\advance\resultno by 1
\bigskip\noindent{\bf Theorem \the\sectno.\the\resultno. }\sl}
\def\defn{ \global\advance\resultno by 1
\bigskip\noindent{\it Definition \the\sectno.\the\resultno. }\slrm}
\def\endthm{\rm\bigskip}
\def\thmend{\rm\bigskip}

\def\endprop{\rm\bigskip}

\def\pf{\rm\medskip\noindent{\it Proof. }}
\def\endpf{\qed\hfil\bigskip}
\def\pfend{\qed\hfil\bigskip}


\def\qed{\hbox{\hskip 1pt\vrule width4pt height 6pt depth1.5pt \hskip
1pt}}

\def\sqr#1#2{{\vcenter{\vbox{\hrule height.#2pt
\hbox{\vrule width.#2pt height#1pt \kern#1pt
\vrule width.2pt}
\hrule height.2pt}}}}



\def\formula{\global\advance\resultno by 1
\eqno{(\the\sectno.\the\resultno)}}
\def\formulano{\global\advance\resultno by 1
(\the\sectno.\the\resultno)}
\def\tableno{\global\advance\resultno by 1
\the\sectno.\the\resultno. }
\def\lformula{\global\advance\resultno by 1
\leqno(\the\sectno.\the\resultno)}


\def\monthname {\ifcase\month\or January\or February\or March\or
April\or
May\or June\or
July\or August\or September\or October\or November\or December\fi}

\newcount\mins  \newcount\hours  \hours=\time \mins=\time
\def\now{\divide\hours by60 \multiply\hours by60 \advance\mins by-\hours

\divide\hours by60 
\ifnum\hours>12 \advance\hours by-12
  \number\hours:\ifnum\mins<10 0\fi\number\mins\ P.M.\else
  \number\hours:\ifnum\mins<10 0\fi\number\mins\ A.M.\fi}



\nopagenumbers
\def\runningtitle{\smallcaps rook monoid algebras and hecke algebras}
\headline={\ifnum\pageno>1\eoheadline\else\firstheadline\fi}
\def\names{\smallcaps tom halverson and arun ram}
\def\firstheadline{}
\def\eoheadline{\ifodd\pageno\oddheadline\else\evenheadline\fi}
\def\oddheadline{\tenrm\hfil\runningtitle\hfil\folio}
\def\evenheadline{\tenrm \folio\hfil{\names}\hfil}

\vphantom{$ $}  
\vskip.75truein

\centerline{\titlefont q-rook monoid algebras, Hecke algebras, and
Schur-Weyl duality}
\bigskip
\centerline{\rm Tom Halverson${}^\ast$}
\centerline{Department of Mathematics and Computer Science}
\centerline{Macalester College}
\centerline{St.\ Paul, MN 55105}
\centerline{{\tt halverson@macalester.edu}}
\medskip
\centerline{\rm and}
\medskip
\centerline{\rm Arun Ram${}^{\ast\ast}$}
\centerline{Department of Mathematics}
\centerline{University of Wisconsin, Madison}
\centerline{Madison, WI 53706}
\centerline{{\tt ram@math.wisc.edu}}

\footnote{}{\tinyrm ${}^\ast$ Research supported in part by National
Science Foundation grant DMS-0100975.}

\footnote{}{\tinyrm \noindent
${}^{\ast\ast}$ Research supported in part by National
Science Foundation grant DMS-9971099, the National Security Agency
and EPSRC grant GR K99015.}

\bigskip



\bigskip
\centerline{\sl  In memory of Sergei Kerov 1946-2000}

\bigskip

{\noindent
\narrower{\narrower{\sl
When we were at the
beginnings of our careers Sergei's support
helped us to believe in our work.  He
generously encouraged us to publish our results on Brauer and
Birman-Murakami-Wenzl
algebras, results which had in part, or possibly
in total, been obtained earlier by Sergei himself.  He remains a great
inspiration
for us, both mathematically and in our memory of his kindness, modesty,
generosity, and encouragement to the younger generation.}\par}\par}
\smallbreak\par

\bigskip

\section 0. Introduction

The rook monoid $R_k$ is the monoid of $k \times k$ matrices with
entries from $\{0,1\}$ and at most one nonzero entry in each row and
column.  Recently, the representation theory of its ``Iwahori-Hecke"
algebra $R_k(q)$, called the  $q$-rook monoid algebra, has been
analyzed. In particular, a Schur-Weyl  type duality on tensor space
was found for the $q$-rook monoid algebra and its irreducible
representations were given explicit combinatorial constructions.  In
this paper we show that, in fact, the $q$-rook monoid algebra is a
quotient of the affine Hecke algebra of type A.  With this knowledge
in hand, we show that the recent results on the $q$-rook monoid
algebras actually come from known results about the affine Hecke
algebra.  In particular

\smallskip\noindent
\itemitem{(a)}  The recent combinatorial construction of the
irreducible representations of $R_k(q)$ by Halverson [Ha] turns out to
be a special case of the construction of irreducible calibrated
representations of  affine Hecke algebras of Cherednik [Ch]
(see also Ram [Ra]), the construction of irreducible representations of
cyclotomic Hecke algebras by Ariki and Koike [AK], and the construction
of the irreducible representations of Iwahori-Hecke algebras of type
B by Hoefsmit [Ho].
\smallskip\noindent
\itemitem{(b)}  The Schur-Weyl duality for the $q$-rook monoid algebra
discovered by Solomon [So2-4] and studied by Halverson [Ha] turns out
to be a special case of the Schur-Weyl duality for cyclotomic Hecke
algebras given by Sakamoto and Shoji [SS].
\smallskip\noindent
Though these results show that the representation theory of the
$q$-rook monoid algebra is ``just'' a piece of the representation
theory of the affine Hecke algebra, this was not at all obvious at
the outset.  It was only on the analysis of the recent
results in [So4] and [Ha] that the similarity to affine
Hecke algebra theory was noticed.  This observation then led us
to search for and establish a concrete connection between these
algebras.

The $q$-rook monoid algebra was first studied in its $q=1$ version in
the 1950's by Munn [Mu1-2].  Solomon [So1] discovered the general
$q$-version of the algebra as a Hecke algebra (double coset algebra)
for the finite algebraic monoid $M_n(\FF_q)$ of $n\times n$ matrices
over a finite field with $q$ elements, with respect to the ``Borel
subgroup'' $B$ of invertible upper triangular matrices.
Later Solomon [So2] found a Schur-Weyl duality for $R_k(1)$ in which
$R_k(1)$ acts as the centralizer algebra for the action of the general
linear group $GL_n(\CC)$ on $V^{\otimes k}$
where $V=L(\varepsilon_1)\oplus L(0)$ is the direct sum of the
``fundamental'' $n$-dimensional representation and the trivial  module
$L(0)$ for $GL_n(\CC)$.  Then Solomon [So3,4]  gave a presentation of
$R_k(q)$ by generators and relations and defined an action of $R_k(q)$
on tensor space.

Halverson [Ha] found a new presentation of $R_k(q)$ and used it to
show that Solomon's action of $R_k(q)$ on tensor space extends the
Schur-Weyl duality so that $R_k(q)$ is the centralizer of the quantum
general linear group $U_q\fgl(n)$ on $V^{\otimes k}$  where
now $V=L(\varepsilon_1)\oplus L(0)$ is the direct sum of the
``fundamental'' and the trivial  module for $U_q\fgl(n)$.
Halverson also  exploited his new presentation to construct,
combinatorially, all the irreducible representations of
$R_k(q)$ when $R_k(q)$ is semisimple.

The main results of this paper are the following:
\smallskip\noindent
\itemitem{(a)} We find yet another presentation (1.6) of $R_k(q)$ by
generators and relations.
\smallskip\noindent
\itemitem{(b)} Our new presentation shows that
$$R_k(q) = H_k(0,1;q)/I,$$
where $H_k(0,1;q)$ is the Iwahori-Hecke algebra of
type $B_k$ with
parameters specialized to $0$ and $1$, and $I$ is the
ideal generated by the  minimal ideal of $H_2(0,1;q)$
corresponding to the pair of partitions
$\lambda = ((1^2),\emptyset)$.
\smallskip\noindent
\itemitem{(c)}  We show that the irreducible representations of
$R_k(q)$ found in [Ha] come from the constructions of irreducible
representations of $H_k(0,1;q)$.
\smallskip\noindent
\itemitem{(d)}  We use the fact that $R_k(q)$ is a quotient of
$H_k(0,1;q)$ and the fact that the Iwahori-Hecke algebra $H_k(q)$ of
type $A_{k-1}$ is a quotient of $R_k(q)$ to easily determine, in
Corollary 2.21, the values of $q$ for which $R_k(q)$ is semisimple.
These values were first found in [So4] using other methods.
\smallskip\noindent
\itemitem{(e)}  We show that the Schur-Weyl duality between $R_k(q)$
and $U_q\fgl(n)$ comes from the Schur-Weyl
duality of Sakamoto and Shoji [SS] for the cyclotomic Hecke algebras
(Theorem 3.5).
\smallskip\noindent
\itemitem{(f)}  We give a {\it different} Schur-Weyl duality for
algebras $A_k(u_1,u_2;q)=H_k(u_1,u_2;q)/I$, where $u_1,u_2\ne 0$,
$H_k(u_1,u_2;q)$ is the Iwahori-Hecke
algebra of type $B_k$ and $I$ is the ideal generated by the
minimal ideal of $H_2(u_1,u_2;q)$ corresponding to the pair of
partitions $\lambda = ((1^2),\emptyset)$.  This Schur-Weyl duality
comes from the Schur-Weyl duality of Orellana and Ram for the affine
Hecke algebra (Theorem 3.3).

\bigskip\noindent
{\bf Acknowledgements.}  A.\ Ram thanks the Isaac Newton Institute for
the Mathematical Sciences at Cambridge University for support
for a very pleasant residency during which the research in this paper
was completed.

\section 1. Presentations of the $q$-rook monoid algebras

Fix $q \in \CC^\ast$. The {\it $q$-rook monoid algebra} is the algebra
$R_k(q)$ given by generators
$$
P_1, P_2, \ldots, P_k \qquad \hbox{and} \qquad
T_1, T_2, \ldots, T_{k-1} $$
with relations
$$
\matrix{
{\rm (A1)} & T_i^2 = (q-q^{-1})T_i+1,\hfill
& 1\le i\le k-1, \hfill \cr
{\rm (A2)} & T_iT_{i+1}T_i = T_{i+1}T_iT_{i+1}, \hfill
& 1\le i\le k-2, \hfill \cr
{\rm (A3)} & T_iT_j = T_jT_i, \hfill & |i-j|>1,\hfill \cr
{\rm (R1)} & P_i^2=P_i, \hfill & 1\le i\le k, \hfill \cr
{\rm (R2)} & P_iP_j=P_jP_i, \hfill & 1\le i,j\le k, \hfill \cr
{\rm (R3)} & P_iT_j = T_jP_i, \hfill & 1\le i<j\le k, \hfill \cr
{\rm (R4)} & P_iT_j = T_jP_i = qP_i, \hfill & 1\le j<i\le k, \hfill \cr
{\rm (R5)} & P_{i+1} = q P_i T_i^{-1} P_i
=q(P_iT_iP_i - (q-q^{-1}) P_i),
\hfill & 1\le i\le k-1.\hfill \cr}\hskip1truein\formula
$$

The algebra $R_k(q)$ was introduced by Solomon [So1] as an analogue of
the Iwahori-Hecke algebra for the finite algebraic monoid $M_k(\FF_q)$
of $k\times k$ matrices over a finite field with $q$ elements with
respect to its ``Borel subgroup'' of invertible upper triangular
matrices.
The presentation of $R_k(q)$ given above is due to Halverson [Ha].

When $q=1$, $R_k(q)$ specializes to the algebra of the rook monoid
$R_k$ that consists of $k\times k$ matrices with entries from
$\{0,1\}$ and {\it at most}
one nonzero entry in each row and column.  These correspond
with the possible placements of nonattacking rooks on an
$k\times k$ chessboard.
In this specialization, $T_i$ becomes the matrix obtained by
switching rows $i$ and $i+1$ in
the identity matrix $I$ and, for $1 \le i \le k-1$, $P_i$ becomes the
matrix $E_{i+1,i+1} + E_{i+1,i+2} +\cdots + E_{k,k}$,
where $E_{i,j}$ is the matrix
with a 1 in position $(i,j)$ and zeros elsewhere.
The generator $P_k$ specializes to
the 0 matrix (which is not the 0 element in the monoid algebra).

\remark The definition of $R_k(q)$ in [So1,3,4] and
[Ha] uses generators $\tilde T_i$ in place of $T_i$.
These generators satisfy
$\tilde T_i^2 = (q-1) \tilde T_i + q$ in place of (A1).
In our presentation, if we let
$\tilde T_i = q T_i$, then
$\tilde T_i^2 = q^2((q-q^{-1})T_i+1)=(q^2-1)qT_i+q^2
=(q^2-1)\tilde T_i+q^2$,
which shows that our algebra is the same except with
parameter $q^2$ instead of $q$.
\thmend

Define
$$
X_i =
T_{i-1} T_{i-2} \cdots T_1(1-P_1)T_1 T_2 \cdots T_{i-1},
\qquad 1 \le i \le k, \formula
$$
so that $X_{i+1} =T_i X_i T_i$.

\lemma In $R_k(q)$ we have the following relations
\smallskip\noindent
\item{{\rm (a)}}  $P_i P_j = P_jP_i=P_j$, for $i\le j$.
\smallskip\noindent
\item{{\rm (b)}}  $P_1 X_2 = P_1 - P_2$.

\pf (a)  If $i=j$ this is (R1).  If $i<j$ then, by (R5) and
induction,
$$P_iP_j = P_i(P_{j-1}T_{j-1}P_{j-1}-(q-q^{-1})P_{j-1})
=P_{j-1}T_{j-1}P_{j-1}-(q-q^{-1})P_{j-1}=P_j.$$
\smallskip\noindent
(b) We use relations (A1), (R4), and (R5) to get
$$
\eqalign{
 P_1 X_2 &=  P_1 T_1 (1 - P_1) T_1
=  P_1 T_1^2 - P_1 T_1 P_1 T_1\cr
&= (q-q^{-1}) P_1 T_1 + P_1 - P_1 T_1 P_1 T_1\cr
&= (q-q^{-1}) P_1 T_1 + P_1 - q^{-1} P_2 T_1 - (q-q^{-1}) P_1 T_1
\cr
&=  P_1 - q^{-1} P_2 T_1
= P_1 - P_2. \hskip.5truein\qed\cr
}
$$

\prop The $q$-rook monoid algebra $R_k(q)$ is generated by
the elements $X_1, T_1, \ldots, T_{k-1}$, and these
elements satisfy the relations (A1), (A2), (A3), and
\smallskip\noindent
\itemitem{{\rm (B1)}} $X_1 T_j = T_j X_1$, for $2\le j\le k$,
\smallskip\noindent
\itemitem{{\rm (B2)}} $X_1^2 = X_1$,
\smallskip\noindent
\itemitem{{\rm (B3)}} $X_1 T_1 X_1 T_1 = T_1 X_1 T_1 X_1$,
\smallskip\noindent
\itemitem{{\rm (B4)}} $(1-X_1)(T_1-q)(1-X_1)(1-X_2) = 0$,
where $X_2 = T_1 X_1 T_1$.
\endprop

\pf By (R5), $R_n(q)$ is generated by
$X_1 = 1-P_1, T_1, \ldots, T_{k-1}$.
\smallskip\noindent
(B1) By (R3), $X_1 T_j = (1-P_1) T_j = T_j(1-P_1) = T_j X_1$.
\smallskip\noindent
(B2) By (R1), $X_1^2 = (1-P_1)^2 = 1 - 2 P_1 + P_1^2 = 1 - P_1 = X_1$.
\smallskip\noindent
(B3) Using (R4) and (R5),
$$
q(T_1 P_1 T_1 P_1 - (q-q^{-1})T_1 P_1)
= T_1 P_2 = P_2 T_1
= q(P_1 T_1 P_1 T_1 - (q-q^{-1})P_1T_1),
$$
and so
$$
P_1 T_1 P_1 T_1 =T_1 P_1 T_1 P_1 + (q-q^{-1}) ( P_1 T_1 - T_1 P_1 ).
\eqno(\ast)
$$
Now, using $(\ast)$ and (A1),
$$
\eqalign{
X_1 T_1 X_1 T_1
&= (1-P_1) T_1 (1 - P_1) T_1 \cr
&= T_1^2 - P_1 T_1^2 - T_1 P_1 T_1 + P_1 T_1 P_1 T_1\cr
&= T_1^2 - (q-q^{-1}) P_1 T_1 - P_1 - T_1 P_1 T_1
+ T_1 P_1 T_1 P_1 + (q-q^{-1})( P_1 T_1 - T_1 P_1)\cr
&= T_1^2 - P_1 - T_1 P_1 T_1 + T_1 P_1 T_1 P_1 - (q-q^{-1})T_1 P_1\cr
&= T_1^2 - T_1^2 P_1 - T_1 P_1 T_1 + T_1 P_1 T_1 P_1 \cr
&= T_1 (1-P_1) T_1 (1-P_1) \cr
&= T_1 X_1 T_1 X_1. \cr
}
$$
Finally, to show (B4), we use Lemma 1.4(b),
$$
\eqalign{
(1-X_1)(T_1-q)(1-X_1)(1-X_2)
& =  P_1(T_1-q)P_1(1-X_2)\cr
& =  P_1T_1P_1 - qP_1 - P_1 T_1 P_1 X_2 + q P_1 X_2\cr
& =  P_1T_1P_1 - qP_1 - P_1 T_1 (P_1-P_2) + q (P_1 - P_2)\cr
& =  -q P_2 + P_1 T_1 P_2\cr
& =  -q P_2 + q P_1 P_2 \qquad\hbox{by (R4)}\cr
& =  -q P_2 + q  P_2 \qquad\hbox{by Lemma 1.4(a)}\cr
& =  0.\qquad\hbox{\qed}\cr
}
$$

Define a new algebra $A_k(q)$ by generators
$$
X_1 \qquad\hbox{and}\qquad T_1, T_2, \ldots, T_{k-1}
$$
and relations
$$
\matrix{
{\rm (A1)} & T_i^2 = (q-q^{-1})T_i+1,\hfill
& 1\le i\le k-1, \hfill \cr
{\rm (A2)} & T_iT_{i+1}T_i = T_{i+1}T_iT_{i+1}, \hfill
& 1\le i\le k-2, \hfill \cr
{\rm (A3)} & T_iT_j = T_jT_i, \hfill & |i-j|>1,\hfill \cr
{\rm (B1)} & X_1 T_j = T_j X_1, \hfill & 2\le j\le k, \hfill \cr
{\rm (B2)} & X_1^2 = X_1, \hfill &  \hfill \cr
{\rm (B3)} & X_1 T_1 X_1 T_1 = T_1 X_1 T_1 X_1, \hfill &\hfill \cr
{\rm (B4)} & (1-X_1)(T_1-q)(1-X_1)(1-X_2) = 0, \hfill
& \hbox{where $X_2 = T_1 X_1 T_1$.}\hfill \cr}\hskip1truein\formula
$$

\medskip\noindent
We will show that $R_k(q) \cong A_k(q)$.  Define
$$
P_1 = (1-X_1) \qquad\hbox{and}\qquad
P_{i+1} = q(P_i T_i P_i - (q-q^{-1}) P_i),
\quad 1 \le i \le k-1. \formula
$$

\lemma Relations (B1)-(B4) are equivalent, respectively, to
\itemitem{{\rm (B1${}'$)}} $P_1 T_j = T_j P_1$, for $2\le j\le n$,
\smallskip\noindent
\itemitem{{\rm (B2${}'$)}} $P_1^2 = P_1$,
\smallskip\noindent
\itemitem{{\rm (B3${}'$)}} $P_2 T_1 = T_1 P_2$,
\smallskip\noindent
\itemitem{{\rm (B4${}'$)}} $P_2^2 = P_2$.
\pf Subtracting $T_j$ from each side of
$$T_j - P_1 T_j = (1-P_1) T_j = X_1 T_j =  T_jX_1
= T_j (1-P_1)  = T_j - T_j P_1,$$
shows that (B1) is equivalent to (B1${}'$).
Relations (B2) and (B2${}'$) are equivalent since
$$1-P_1 =  X_1 = X_1^2 = (1-P_1)^2 = 1 - 2 P_1 + P_1^2.$$
Since
$$
\eqalign{
X_1 T_1 X_1 T_1
&= (1-P_1) T_1 (1-P_1) T_1
= (1-P_1)T_1^2 -  T_1 P_1 T_1 + P_1 T_1 P_1 T_1\cr
&= (1-P_1)((q-q^{-1})T_1 + 1) -T_1 P_1 T_1
+(q^{-1}P_2+(q-q^{-1})P_1)T_1\quad\hbox{(by (A1) and (1.7))} \cr
&= (q-q^{-1})T_1+1 -(q-q^{-1})P_1T_1-P_1
- T_1 P_1 T_1 + q^{-1}P_2T_1 + (q-q^{-1}) P_1 T_1\cr
&= (q-q^{-1})T_1 + 1 - P_1 - T_1 P_1 T_1 + q^{-1}P_2T_1\cr
}
$$
is equal to
$$
\eqalign{
T_1 X_1 T_1 X_1
&= T_1 (1-P_1) T_1 (1-P_1)
= T_1^2(1-P_1) - T_1P_1T_1 + T_1P_1T_1P_1\cr
&=((q-q^{-1})T_1+1)(1-P_1) -T_1P_1T_1 +T_1(q^{-1}P_2+(q-q^{-1})P_1)
\quad\hbox{(by (A1) and (1.7))} \cr
&=(q-q^{-1})T_1+1 -(q-q^{-1})T_1P_1-P_1 -T_1P_1T_1 +q^{-1}T_1P_2
+ (q-q^{-1})T_1P_1 \cr
&= (q-q^{-1})T_1+1 -P_1 -T_1P_1T_1 +q^{-1}T_1P_2,\cr
}
$$
(B3) is equivalent to (B3$'$).

\smallskip
Expanding
$$
\eqalign{
(T_1 - q)(1-X_2) &= (T_1 - q) (1 - T_1 X_1 T_1)
= (T_1 - q) (1 - T_1 (1-P_1)  T_1)  \cr
&= (T_1 - q) ( 1 - (q-q^{-1}) T_1 - 1 + T_1P_1T_1)
\qquad\hbox{(by (A1))} \cr
&= (T_1 - q)T_1(-(q-q^{-1}) + P_1T_1)  \cr
&= ((q-q^{-1})T_1 + 1) - qT_1)(q^{-1}-q + P_1T_1)
\qquad\hbox{(by (A1))} \cr
&= (1-q^{-1}T_1)(q^{-1}-q + P_1T_1)  \cr
&= q^{-1}-q + P_1T_1 -q^{-2}T_1+T_1 -q^{-1}T_1P_1T_1  \cr
&= -(q-q^{-1}) + q^{-1}(q-q^{-1})T_1 + P_1T_1- q^{-1}T_1P_1T_1, \cr
}
$$
gives
$$
\eqalign{
(1-&X_1)(T_1 - q)(1-X_2)(1-X_1)
= P_1(T_1 - q)(1-X_2)P_1 \cr
&= -(q-q^{-1})P_1 + q^{-1}(q-q^{-1})P_1T_1P_1 + P_1T_1P_1
- q^{-1}P_1T_1P_1T_1P_1
\quad\hbox{(by (B2${}'$))}\cr
&= -(q-q^{-1})P_1 + (2-q^{-2})P_1T_1P_1 - q^{-1}(P_1T_1P_1)^2
\quad\hbox{(by (B2${}'$))}\cr
&= -(q-q^{-1})P_1 + (2-q^{-2})(q^{-1}P_2 + (q-q^{-1})P_1)
-q^{-1}(q^{-1}P_2 + (q-q^{-1})P_1)^2
\quad\hbox{(by (1.7))}\cr
&= (q-q^{-1})(-1+2-q^{-2})P_1  + q^{-1}(2-q^{-2})P_2 \cr
&\qquad\qquad\qquad\qquad
-q^{-1}(q^{-2}P_2^2 + 2q^{-1}(q-q^{-1}) P_2 +(q-q^{-1})^2P_1) \cr
&= (q-q^{-1})^2(q^{-1}-q^{-1})P_1
+q^{-1}(q^{-2}P_2 - q^{-2} P_2^2) \cr
&=q^{-3}(P_2 -  P_2^2). \cr
}
$$
and so $(1-X_1)(T_1 - q)(1-X_2)(1-X_1)=0$ if and only if $P_2^2 = P_2$.
Thus (B4) is equivalent to (B4${}'$).
\pfend

\prop  The algebra $A_k(q)$ is generated by
$T_1, \ldots, T_{k-1}, P_1, \ldots, P_k$, and these elements satisfy
the relations (A1), (A2), (A3), and
\itemitem{(E1)} $P_iT_j = T_jP_i$, for $1\le i<j\le k$,
\smallskip\noindent
\itemitem{(E2)} $P_iP_j=P_jP_i = P_i$, for all $1\le j < i \le k$,
\smallskip\noindent
\itemitem{(E3)} $P_i^2=P_i$, for $1\le i\le k$,
\smallskip\noindent
\itemitem{(E4)} $P_iT_j = T_jP_i = qP_i$, for $1\le j<i\le k$,
\smallskip\noindent
\itemitem{(E5)} $P_{i+1} = q P_i T_i^{-1} P_i = q (P_iT_iP_i -
(q-q^{-1})
P_i),$ for $1\le i\le k-1$.

\pf Since $X_1 = 1 - P_1$, the elements
$T_1, \ldots, T_{n-1}, P_1$ generate $A_k(q)$.
Relation (E5) is the definition of $P_{i+1}$.

We prove (E1) by induction on $i$.  The case $i=1$ is (B1${}'$).
Assume that $j > i+1 > 1$, then $T_j$ commutes with $P_i$ by
induction and $T_j$ commutes with $T_i$ by (A3), so
$$
T_j P_{i+1} = q^{-1}T_jP_i T_j P_i - (q-q^{-1}) T_j P_i
= q^{-1}P_i T_j P_i T_j- (q-q^{-1}) P_i T_j =  P_{i+1} T_j,
$$
proving (E1).

We now prove (E2)-(E4) collectively by induction on $i$.
The case $i=1$ for (E2) follows from (B2${}'$) since
$$
P_2 P_1 = q(P_1 T_1 P_1 - (q-q^{-1}) P_1)P_1
= q(P_1 T_1 P_1 - (q-q^{-1}) P_1) = P_2.
\eqno(\ast)
$$
The relation $P_1 P_2 = P_2$ is similar.
The first two cases of (E3) are (B2${}'$) and (B4${}'$).
The $i=1$ case of (E4) follows from (B4${}'$) since
$$
\eqalign{
P_2^2 &= P_2q(P_1T_1P_1 - (q-q^{-1})P_1) \cr
&= q(P_2T_1P_1 - (q-q^{-1})P_2) \qquad\hbox{(by ($\ast$))}\cr
&= q(T_1P_2P_1 - (q-q^{-1})P_2) \qquad\hbox{(by (B3${}'$))}\cr
&= q(T_1P_2 - (q-q^{-1})P_2) \qquad\hbox{(by ($\ast$))}.\cr
}
$$
Since $P_2^2 = P_2$, we have $T_1 P_2 = qP_2$.
The case $P_2 T_1 = qP_2$ is similar.

Now fix $i > 1$ and assume the following relations,
\itemitem{(E2${}^\ast$)} $P_iP_j=P_jP_i = P_i$, for all $1\le j < i$,
\smallskip\noindent
\itemitem{(E3${}^\ast$)} $P_i^2=P_i$,
\smallskip\noindent
\itemitem{(E4${}^\ast$)} $P_iT_j = T_jP_i = qP_i$, for $1\le j<i$,
\smallskip\noindent
We show each of these relations for $i+1$.

For (E2) we use  (E3${}^\ast$) to get
$$
P_{i+1} P_i = q(P_i T_i P_i^2 - (q-q^{-1}) P_i^2)
= q(P_i T_i P_i - (q-q^{-1}) P_i) = P_{i+1},
$$
and when $j < i$, we use (E2${}^\ast$) to get
$$
P_{i+1} P_j = q(P_i T_i P_i P_j - (q-q^{-1}) P_i P_j)
= q(P_i T_i P_i - (q-q^{-1}) P_i) = P_{i+1}.
$$
The relations $P_i P_{i+1} = P_{i+1}$ and $P_j P_{i+1} = P_{i+1}$ are
similar, and so (E2) is established.

To establish (E3), let $i>2$
(note that we have established $i=1,2$),
$$
\eqalign{
&P_{i+1}^2
= q^2(P_i T_i P_i - (q-q^{-1})P_i)^2 \cr
&= q^2(P_i T_i P_i T_i P_i - 2(q-q^{-1})P_iT_iP_i
+(q-q^{-1})^2 P_i) \quad\hbox{(by (E3${}^\ast$))}\cr
&= q^2(P_i T_i qP_{i-1} T_{i-1} P_{i-1} T_i P_i
- q(q-q^{-1})P_i T_i P_{i-1} T_i P_i \cr
&\qquad\qquad\qquad\qquad
- 2(q-q^{-1}) P_iT_iP_i + (q-q^{-1})^2 P_i) \quad\hbox{(by (E5))}\cr
&= q^2(qP_iP_{i-1}T_iT_{i-1}T_iP_{i-1}P_i
-q(q-q^{-1})P_iP_{i-1}T_i^2P_i \cr
&\qquad\qquad\qquad\qquad
-2(q-q^{-1})P_iT_iP_i +(q-q^{-1})^2 P_i) \quad\hbox{(by (E1))}\cr
&= q^2(qP_iT_iT_{i-1}T_iP_i -q(q-q^{-1})P_i T_i^2 P_i
-2(q-q^{-1})P_iT_iP_i +(q-q^{-1})^2 P_i) \quad\hbox{(by (E2))}\cr
&=q^2(qP_iT_{i-1}T_iT_{i-1}P_i
-q(q-q^{-1})P_i T_i^2 P_i
-2(q-q^{-1})P_iT_iP_i +(q-q^{-1})^2 P_i) \quad\hbox{(by (A2))}\cr
&= q^2(qP_iT_{i-1}T_iT_{i-1}P_i
-q(q-q^{-1})^2P_i T_i P_i - q(q-q^{-1})P_i \cr
&\qquad\qquad\qquad\qquad
-2(q-q^{-1})P_iT_iP_i +(q-q^{-1})^2 P_i \quad\hbox{(by (A1))}\cr
&= q^2(q^3 P_iT_iP_i
-q(q-q^{-1})^2P_i T_i P_i - q(q-q^{-1})P_i \cr
&\qquad\qquad\qquad\qquad
-2(q-q^{-1}) P_iT_iP_i +(q-q^{-1})^2 P_i)
\quad\hbox{(by (E4${}^\ast$))}\cr
&= q^2((q^3-q(q-q^{-1})^2-2(q-q^{-1})) P_iT_iP_i
+ (q-q^{-1})(-q+q-q^{-1}) P_i) \cr
&= q(P_iT_iP_i - (q-q^{-1}) P_i) \cr
&= P_{i+1}. \cr
 }
$$

Finally, we prove (E4).  First let $j < i$.  Then,
by (E4${}^\ast$),
$$
P_{i+1} T_j =  q(P_iT_iP_i T_j - (q-q^{-1})P_iT_j)
= q(qP_iT_iP_i - (q-q^{-1}) q P_i) = q P_{i+1},
$$
and $T_jP_{i+1}  = q P_{i+1}$ is similar.  Now we consider the case
where $j = i$,
$$
\eqalign{
P_{i+1} T_i
&= q(P_i T_i P_i T_i - (q-q^{-1}) P_i T_i)\quad\hbox{(by (E5))}\cr
&= q(P_iT_iqP_{i-1}T_{i-1}P_{i-1}T_i -q(q-q^{-1})P_iT_iP_{i-1}T_i
-(q-q^{-1})P_iT_i)\quad\hbox{(by (E5))}\cr
&= q(qP_iP_{i-1}T_iT_{i-1}T_iP_{i-1}
-q(q-q^{-1})P_iP_{i-1}T_i^2
-(q-q^{-1})P_iT_i)\quad\hbox{(by (E1))}\cr
&=q(qP_iT_iT_{i-1}T_iP_{i-1}
-q(q-q^{-1})P_iT_i^2 -(q-q^{-1})P_iT_i)\quad\hbox{(by (E2))}\cr
&= q(qP_iT_{i-1}T_iT_{i-1}P_{i-1}
-q(q-q^{-1})P_iT_i^2 -(q-q^{-1})P_iT_i)
\quad\hbox{(by (A2))}\cr
&= q(q^2P_iT_iT_{i-1}P_{i-1}  -q(q-q^{-1})P_iT_i^2 -(q-q^{-1}) P_iT_i)
\quad\hbox{(by (E4${}^\ast$))}\cr
&= q(q^2P_iT_iT_{i-1}P_{i-1} -q(q-q^{-1})^2P_iT_i -q(q-q^{-1})P_i
-(q-q^{-1}) P_iT_i)\quad\hbox{(by (A1))}\cr
&= q(q^2P_iT_iT_{i-1}P_{i-1} -q^2(q-q^{-1})P_iT_i
-q(q-q^{-1}) P_i) \cr
&= q(q^2P_iP_{i-1}T_iT_{i-1} P_{i-1}  -q^2(q-q^{-1})P_iT_i
-q(q-q^{-1})P_i)\quad\hbox{(by (E2))} \cr
&= q(q^2P_iT_iP_{i-1}T_{i-1} P_{i-1}  -q^2(q-q^{-1})P_iT_i
-q(q-q^{-1})P_i)\quad\hbox{(by (E1))} \cr
&= q(qP_iT_iP_i +q^2(q-q^{-1})P_iT_iP_{i-1} -q^2(q-q^{-1})P_iT_i
-q(q-q^{-1})P_i)\quad\hbox{(by (E5))} \cr
&= q(P_{i+1} +q^2(q-q^{-1})P_iT_iP_{i-1} -q^2(q-q^{-1})P_iT_i)
\quad\hbox{(by (E5))} \cr
&= q(P_{i+1} +q^2(q-q^{-1})P_iP_{i-1}T_i -q^2(q-q^{-1})P_iT_i)
\quad\hbox{(by (E1))} \cr
&= q(P_{i+1} +q^2(q-q^{-1})P_iT_i -q^2(q-q^{-1})P_iT_i)
\quad\hbox{(by (E2))} \cr
&= qP_{i+1}. \cr
}
$$
The case $T_i P_{i+1} = qP_{i+1}$ is similar.  \pfend

Propositions 1.5 and 1.9 give the following theorem.

\thm $R_k(q) \cong A_k(q)$ and thus (1.6) is a new presentation of
$R_k(q)$.
\endthm

\section 2. Hecke algebras

\subsection The affine Hecke algebra $\tilde H_k$

Fix $q\in \CC^*$.  The {\it affine Hecke algebra} $\tilde H_k$ is the
algebra given by generators
$$X_1, \ldots, X_k
\qquad\hbox{and}\qquad
T_1,\ldots, T_{k-1}$$
with relations
$$
\matrix{
{\rm (1)} & T_iT_j=T_jT_i, \hfill & |i-j|>1, \hfill\cr
{\rm (2)} & T_iT_{i+1}T_i=T_{i+1}T_iT_{i+1}, \hfill& 1\le i\le
k-2,\hfill\cr
{\rm (3)} & T_i^2 = (q-q^{-1})T_i+1, \hfill& 1\le i\le k,\hfill\cr
{\rm (4)} & X_iX_j = X_jX_i,\hfill & 1\le i,j\le k,\hfill\cr
{\rm (5)} & X_iT_i = T_iX_{i+1}+(q-q^{-1})X_i,\hfill & 1\le i\le
k-1.\hfill\cr}
\hskip2truein
$$

\smallskip\noindent
It follows from relations (3) and (5) that
$$X_i = T_{i-1}\cdots T_2T_1X_1T_1T_2\cdots T_{i-1},
\qquad\hbox{for $1\le i\le k$,}$$
and from (4) that
\smallskip\noindent
$$
\matrix{
{\rm (6)} & X_1T_1X_1T_1 = T_1X_1T_1X_1. \hfill & \hskip1.3truein
\hfill\cr}
\hskip2truein
$$
\smallskip\noindent
In fact, $\tilde H_k$ can be presented as the algebra generated by
$X_1$ and $T_1,\ldots, T_{k-1}$ with relations (1-3) and (6).

Let
$$[k]! =[1][2]\cdots [k],
\qquad\hbox{where}\qquad
[i] = 1+q^2+\cdots+q^{2(i-1)}.$$
When $[k]!\ne 0$ a large class of irreducible representations
of the affine Hecke algebra $\tilde H_k(q)$,
the integrally calibrated irreducible representations,
have a simple combinatorial
construction.  An $\tilde H_k$-module $M$ is {\it integrally
calibrated\/} if $M$ has a basis of simultaneous eigenvectors for
$X_1,X_2,\ldots, X_k$ for which the eigenvalues are all of the form
$q^j$
with $j\in \ZZ$.
The construction of these $\tilde H_k$-modules is originally due to
Cherednik [Ch] (see [Ra] for greater detail) and is a generalization
of the classical seminormal construction of the irreducible
representations
of the symmetric group by A. Young. Young's construction
had been generalized to Iwahori-Hecke algebras of classical type
(see Theorem 2.8 below) by Hoefsmit [Ho] in 1974.

To describe the construction we shall use the notations
of [Mac] for partitions so that
a partition is identified with a collection of boxes in
a corner, $\ell(\lambda)$ is the number of rows of $\lambda$, and
$|\lambda|$ is the number of boxes in $\lambda$. For example,
the partition
$$\lambda=(5,5,3,1,1)=
{\beginpicture
\setcoordinatesystem units <0.25cm,0.25cm>     
\setplotarea x from 0 to 4, y from -1 to 1    
\linethickness=0.05pt      
\putrectangle corners at 0 0  and 1 1
\putrectangle corners at 2 0  and 1 1
\putrectangle corners at 3 0  and 2 1
\putrectangle corners at 0 2  and 1 1
\putrectangle corners at 2 2  and 1 1
\putrectangle corners at 3 2  and 2 1
\putrectangle corners at 4 2  and 3 1
\putrectangle corners at 4 2  and 5 1

\putrectangle corners at 0 2  and 1 3
\putrectangle corners at 2 2  and 1 3
\putrectangle corners at 3 2  and 2 3
\putrectangle corners at 4 2  and 3 3
\putrectangle corners at 4 2  and 5 3

\putrectangle corners at 0 -1  and 1 0
\putrectangle corners at 0 -2  and 1 -1
\endpicture}
\qquad\hbox{has $\ell(\lambda) = 5$ and $|\lambda|=15$.}$$

If $\lambda$ is a partition that is obtained from $\mu$ by adding $k$
boxes let $\lambda/\mu$ be the {\it skew shape} consisting of those
boxes of
$\lambda$ that are not in $\mu$.
A {\it standard tableau} of shape $\lambda/\mu$ is a filling of the
boxes
of $\lambda/\mu$ with $1,2,\ldots, k$ such that
\smallskip\noindent
\itemitem{(a)} the entries in the rows increase left to right, and
\smallskip\noindent
\itemitem{(b)} the entries in the columns increase top to bottom.
\smallskip\noindent
If $b$ is a box in $\lambda/\mu$ define
$$CT(b) = q^{2(c-r)},
\qquad\hbox{if $b$ is in position $(r,c)$ of $\lambda$.}
\formula
$$

\thm ([Ch], see also [Ra, Theorem 4.1], and [OR, Theorem 6.20a])
Assume that $[k]!\ne 0$.  Then the
calibrated irreducible representations $H^{\lambda/\mu}$ of the
affine Hecke algebra $\tilde H_k$ are indexed by skew shapes
and can be given explicitly as the vector space
$$H^{\lambda/\mu} = \hbox{$\CC$-span}\{ v_L\ |\ \hbox{$L$ is a standard
tableau of shape $\lambda/\mu$}\}$$
(so that the symbols $v_L$ form a basis of $H^{\lambda/\mu}$) with
$\tilde H_k$-action given by
$$\eqalign{
X_iv_L &= CT(L(i)) v_L, \cr
T_i v_L &= \left({CT(L(i+1))(q-q^{-1})\over CT(L(i+1))-CT(L(i))}\right)
v_L
+ \left(q^{-1}+{CT(L(i+1))(q-q^{-1})\over CT(L(i+1))-CT(L(i))}\right)
v_{s_iL},
\cr}
$$
where
\smallskip\noindent
\itemitem{} $s_iL$ is the same as $L$ except $i$
and $i+1$ are switched, and
\smallskip\noindent
\itemitem{} $v_{s_iL} =0$, if $s_iL$ is not a standard tableau.
\endthm

\remark
In [OR] it is explained how the basis $v_L$ of $H^{\lambda/\mu}$
and the action of $\tilde H_k$ in Theorem 2.2 can be derived in
a natural way from the general mechanism of quantum groups
(${\cal R}$-matrices, quantum Casimirs,
the tensor product rule in (3.1))
and a Schur-Weyl duality theorem (Theorem 3.3 below) for the affine
Hecke algebra.
\endthm

\subsection The cyclotomic Hecke algebra $H_k(u_1,\ldots,u_r;q)$

Let $u_1, \ldots, u_r \in \CC$ and $q \in \CC^\ast$.  The
{\it cyclotomic Hecke algebra\/} $H_k(u_1,\ldots,u_r;q)$
is the quotient of the affine Hecke algebra $\tilde H_k$
by the ideal generated by the relation
$$(X_1-u_1)(X_1-u_2)\cdots(X_1-u_r)=0.
\formula$$
The algebra $H_k(u_1,\ldots,u_r;q)$ is a deformation of
the group algebra of the complex reflection
group $G(r,1,k) = (\ZZ/r\ZZ) \wr S_k$ and is of dimension
$\dim (H_k(u_1,\ldots,u_r;q))=r^kk!$.
These algebras were introduced by Ariki and Koike [AK].
Ariki has generalized the classical result of Gyoja and Uno [GU]
and given precise conditions for the semisimplicity of
the cyclotomic Hecke algebras.

\thm ([Ar]) The algebra $H_k(u_1,,\ldots, u_r;q)$ is semisimple if
and only if
$$
\hbox{$q^{2d} u_i\ne u_j$ for all $-k<d<k$, $1\le i<j\le r$,}
\qquad\hbox{and}\qquad
\hbox{\ \ $[k]!\ne 0$,}$$
where $[k]! = [1][2]\cdots [k]$ and $[i] = 1+q^2+\cdots +q^{2(i-1)}$.
\endthm
\pf
Let us only explain the conversion between the statement in [Ar]
and the statement here.  This conversion is the same as in Remark 1.2.
If
$\tilde T_i = qT_i$ then
$\tilde T_i^2 = q^2((q-q^{-1})T_i+1)=(q^2-1)qT_i+q^2=(q^2-1)\tilde
T_i+q^2$,
which shows that our algebra is the same as Ariki's except with
parameter $q^2$.
\endpf

Ariki and Koike give a combinatorial construction of
the irreducible representations of the cyclotomic Hecke
algebra $H_k(u_1,\ldots, u_r;q)$
when it is semisimple. Define
$$\hat H_k^{(r)} = \{\hbox{$r$-tuples
$\lambda = (\lambda^{(1)},\ldots,\lambda^{(r)})$ of
partitions with $k$ boxes total}\}.
\formula
$$
Let $\lambda\in \hat H_k^{(r)}$.  A {\it standard tableau of shape
$\lambda$}
is a filling of the boxes of $\lambda$ with $1,2,\ldots, k$ such that
for each $\lambda^{(i)}$, $1\le i\le r$,
\smallskip\noindent
\itemitem{(a)} the entries in the rows are increasing left to right,
and
\smallskip\noindent
\itemitem{(b)} the entries in the columns are increasing top to
bottom.
\smallskip\noindent
Let $L(i)$ denote the the box of $L$ containing $i$, and define
$$
CT(b) = u_iq^{2(c-r)}, \qquad
\hbox{ if the box $b$ is in position $(r,c)$ of $\lambda^{(i)}$.}
\formula
$$

\thm [AK, Theorem 3.7]  If $H_k(u_1,u_2,\ldots, u_r;q)$ is semisimple
its irreducible representations $H^\lambda$,
$\lambda\in \hat H_k^{(r)}$, are given by
$$H^{\lambda} = H^{(\lambda^{(1)},\ldots,\lambda^{(r)})}
= \hbox{$\CC$-span}\{ v_L\ |\
\hbox{$L$ is a standard tableau of shape $\lambda$}\}$$
(so that the symbols $v_L$ form a basis of the vector
space $H^\lambda$) with $H_k(u_1,\ldots,u_r;q)$-action
given by
$$
X_iv_L = CT(L(i)) v_L,
\qquad\hbox{and}\qquad
T_i v_L = (T_i)_{LL} v_L + (q^{-1}+(T_i)_{LL}) v_{s_iL}, $$
where
\smallskip\noindent
\itemitem{} $\displaystyle{
(T_i)_{LL} = \cases{
q, &if $L(i)$ and $L(i+1)$ are in the same row, \cr
& \quad of $\lambda^{(j)}$ for some fixed $j$, \cr
-q^{-1}, &if $L(i)$ and $L(i+1)$ are in the same column, \cr
& \quad of $\lambda^{(j)}$ for some fixed $j$, \cr
\displaystyle{{CT(L(i+1))(q-q^{-1})\over CT(L(i+1))-CT(L(i))}},
&otherwise, \cr} }$
\medskip\noindent
\itemitem{} $s_iL$ is the same as $L$ except $i$
and $i+1$ are switched, and
\medskip\noindent
\itemitem{} $v_{s_iL} =0$, if $s_iL$ is not a standard tableau.
\endthm

\noindent
It is interesting to note that Theorem 2.8
is almost an immediate consequence of Theorem 2.2.

\subsection The Iwahori-Hecke algebra $H_k(u_1,u_2;q)$ of type $B_k$

The {\it Iwahori-Hecke algebra of type $B_k$} is the
cyclotomic Hecke algebra $H_k(u_1,u_2;q)$.
Thus, for $u_1, u_2\in \CC$ and $q\in \CC^*$,
$H_k(u_1,u_2;q)$ is the quotient
of the affine Hecke algebra by the ideal generated by the relation
$$(X_1-u_1)(X_1-u_2)=0. \formula$$
The algebra $H_k(1,-1;q)$ is the group algebra of the Weyl
group of type
$B_k$ (the hyperoctahedral group of signed permutations).

In the case of the Iwahori-Hecke algebra $H_k(u_1,u_2;q)$ of type $B_k$
Theorem 2.8 is due to Hoefsmit [Ho].
When $H_k(u_1,u_2;q)$ is semisimple Hoefsmit's construction of the
irreducible representations of $H_k(u_1,u_2;q)$ implies that,
as $H_{k-1}(u_1,u_2;q)$-modules
$${\rm Res}^{H_k}_{H_{k-1}} H^\lambda = \bigoplus_{\lambda^-}
H^{\lambda^-},\formula$$
where the sum runs over all pairs of partitions $\lambda^-$ which are
obtained
from $\lambda$ by removing a single box, and
$$H^{\lambda^-}= \hbox{$\CC$-span}\left\{ v_L\ \Big|\
\matrix{
\hbox{$L$ is a standard tableau of shape $\lambda$}\hfill\cr
\hbox{and $L^-$ has shape $\lambda^-$}\hfill \cr}
\right\},\formula$$
where $L^-$ is the standard tableau with $k-1$ boxes which is obtained
by removing the entry $k$ from $L$.  The restriction rules (2.10)
can be encoded in the {\it Bratteli diagram} for the sequence of
algebras
$$H_1(u_1,u_2;q)\subseteq H_2(u_1,u_2;q)\subseteq H_3(u_1,u_2;q)
\subseteq \cdots
\formula$$
i.e., the graph which has
$$\hbox{vertices on level $k$ indexed by $\lambda\in \hat H_k^{(2)}$}
\qquad\hbox{and}\qquad
\hbox{edges $\lambda \longleftrightarrow \lambda^-$}$$
if $\lambda^-$ is obtained from $\lambda$ by removing a single box.
The first few rows of the Bratteli diagram for $H_k(u_1,u_2;q)$ are
displayed in Figure 1.

$$ {\beginpicture
\setcoordinatesystem units <0.2cm,0.2cm>     
\setplotarea x from -3 to 33, y from 11 to 30     
\linethickness=0.05pt      

\put{$(\emptyset,\emptyset)$} at 20 30

\put{$(\emptyset, \phantom{t,,} )$} at 12 25
\put{$(\phantom{t,,} ,\emptyset)$} at 28 25
\putrectangle corners at 12.5 24.5 and 13.5 25.5
\putrectangle corners at 26.5 24.5 and 27.5 25.5
\plot 19 28.5 12 26.5 /
\plot 21 28.5 28 26.5 /

\put{$(\emptyset,\phantom{ttt,}  )$} at 4 20
\putrectangle corners at 4 19.5 and 5 20.5
\putrectangle corners at 6 19.5 and 5 20.5

\put{$(\emptyset,\phantom{tt,})$} at 12 20
\putrectangle corners at 12.5 21 and 13.5 20
\putrectangle corners at 12.5 19 and 13.5 20

\put{$(\phantom{tt},\phantom{tt})$} at 20 20
\putrectangle corners at 20.3 19.5 and 21.3 20.5
\putrectangle corners at 18.5 19.5 and 19.5 20.5

\put{$(\phantom{ttt,} ,  \emptyset)$} at 28 20
\putrectangle corners at 26.8 19.5 and 27.8 20.5
\putrectangle corners at 25.8 19.5 and 26.8 20.5

\put{$(\phantom{tt,},  \emptyset)$} at 36 20
\putrectangle corners at 34.5 20 and 35.5 21
\putrectangle corners at 34.5 19 and 35.5 20

\plot  11 23.5  4 21.5 /
\plot  12 23.5  12 21.5 /
\plot  13 23.5  19 21.5 /
\plot  27 23.5  21 21.5 /
\plot  28 23.5  28 21.5 /
\plot  29 23.5  35 21.5 /

\put{$(\! \emptyset,\phantom{ttt,}\, )$} at -4.2 14
\putrectangle corners at -5 13.5 and -4 14.5
\putrectangle corners at -3 13.5 and -2 14.5
\putrectangle corners at -4 13.5 and -3 14.5

\put{$(\! \emptyset,\,\phantom{ttt})$} at 2.2 14
\putrectangle corners at 3 12.5 and 2 13.5
\putrectangle corners at 3 13.5 and 2 14.5
\putrectangle corners at 4 13.5 and 3 14.5

\put{$( \emptyset,\phantom{t}\, )$} at 7.3 14
\putrectangle corners at 8.5 13.5 and 7.5 14.5
\putrectangle corners at 8.5 13.5 and 7.5 12.5
\putrectangle corners at 8.5 11.5 and 7.5 12.5

\put{$( \phantom{t,}\,,\phantom{tt}\,)$} at 12.5 14
\putrectangle corners at 13.5 13.5 and 12.5 14.5
\putrectangle corners at 13.5 13.5 and 14.5 14.5
\putrectangle corners at 11.5 13.5 and 10.5 14.5

\put{$( \phantom{t,},\phantom{t,} )$} at 17.7 14
\putrectangle corners at 19 13.5 and 18 14.5
\putrectangle corners at 19 13.5 and 18 12.5
\putrectangle corners at 17.2 13.5 and 16.2 14.5

\put{$( \phantom{tt,}\, ,\phantom{t} )$} at 22.7 14
\putrectangle corners at 24.5 13.5 and 23.5 14.5
\putrectangle corners at 23 13.5 and 22 14.5
\putrectangle corners at 21 13.5 and 22 14.5

\put{$( \phantom{t,}\, ,\phantom{t} )$} at 27.5 14
\putrectangle corners at 29 13.5 and 28 14.5
\putrectangle corners at 27.3 13.5 and 26.3 12.5
\putrectangle corners at 27.3 13.5 and 26.3 14.5

\put{$(\phantom{ttt,}\,\, ,  \emptyset\!)$} at 32.8 14
\putrectangle corners at 31.5 13.5 and 30.5 14.5
\putrectangle corners at 32.5 13.5 and 31.5 14.5
\putrectangle corners at 33.5 13.5 and 32.5 14.5

\put{$( \phantom{tt,}\, ,\emptyset)$} at 38.5 14
\putrectangle corners at 37.7 13.5 and 36.7 12.5
\putrectangle corners at 37.7 13.5 and 36.7 14.5
\putrectangle corners at 38.7 13.5 and 37.7 14.5

\put{$(\phantom{t,}\,,\emptyset)$} at 44 14
\putrectangle corners at 43.6 12 and 42.6 13
\putrectangle corners at 43.6 13 and 42.6 14
\putrectangle corners at 43.6 14 and 42.6 15

\plot 3 18.5  -3 15.5 /
\plot 4 18.5  2.5 15.5 /
\plot 5 18.5  12 15.5 /
\plot 19 18.5  13 15.5 /
\plot 20 18.5  18 15.5 /
\plot 11 18.5  3.5 15.5 /
\plot 11.5 18.5  8.5 15.5 /
\plot 13.5 18.5  17 15.5 /
\plot 20.5 18.5  22 15.5 /
\plot 21.5 18.5  27 15.5 /
\plot 27 18.5  23 15.5 /
\plot 28 18.5  32 15.5 /
\plot 29 18.5 37 15.5 /
\plot 36 18.5 38 15.5 /
\plot 37 18.5 43 15.5 /
\plot 35 18.5 28 15.5 /
\endpicture}
$$
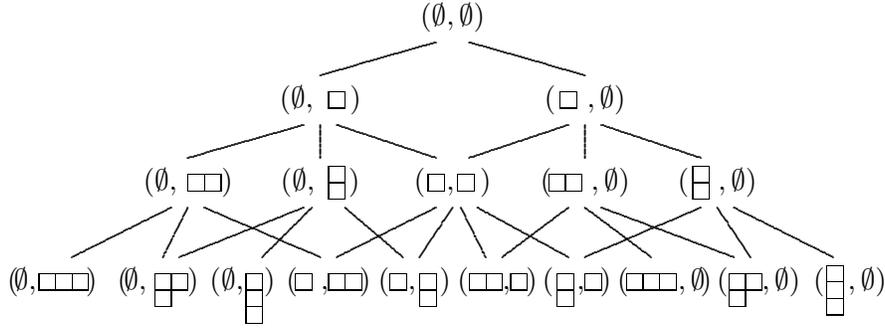
\centerline{Figure 1. Bratteli Diagram for $H_k(u_1,u_2;q)$}

\medskip
If $H_k(u_1,u_2;q)$ is semisimple then
$$H_k(u_1,u_2;q)\cong
\bigoplus_{\lambda\in \hat B_k} M_{d_\lambda}(\CC),
\formula$$
where $d_\lambda$ is the number of standard tableaux $L$ of
shape $\lambda$, and $M_d(\CC)$ is the algebra of
$d\times d$ matrices with entries from $\CC$.

The minimal ideals $I^\lambda$ of $H_k(u_1,u_2;q)$ are in
one-to-one correspondence with the summands in (2.13).
Let $m<k$, let $I^\mu$ be a fixed minimal ideal
of $H_m(u_1,u_2;q)$ and define
$$\hbox{$\langle I^\mu\rangle_k$ is the ideal of $H_k(u_1,u_2;q)$
generated by $I^\mu$}$$
($I^\mu\subseteq H_m(u_1,u_2;q)\subseteq H_k(u_1,u_2;q)$).
The restriction rules (2.10) imply that
$$\langle I^\mu\rangle_k
= \bigoplus_{\lambda\supseteq\mu} I^\lambda, \formula$$
where the sum is over all pairs of partitions
$\lambda = (\lambda^{(1)},\lambda^{(2)})\in H_k^{(2)}$ which are
obtained from
$\mu=(\mu^{(1)},\mu^{(2)})\in \hat H_m^{(2)}$ by adding $(k-m)$ boxes.

\subsection The ideal $I^{((1^2),\emptyset)}$

\lemma Assume that $H_2(u_1,u_2;q)$ is semisimple.
The minimal ideal $I^{((1^2),\emptyset)}$ of $H_2(u_1,u_2;q)$ is
generated
by the element
$$p=\cases{
(X_1-u_2)(X_2-u_2)(X_2-q^2u_1), &if $u_1\ne 0$, \cr
(X_1-u_2)(T_1-q)(X_1-u_2)(X_2-u_2), &if $u_1=0$, \cr
}\qquad\quad\hbox{where $X_2=T_1X_1T_1$.}$$
\pf  Using the construction of the simple $H_2(u_1,u_2;q)$-modules
in Theorem 2.8 it is not tedious to check that, when $u_1\ne 0$,
$$pv_L = \cases{
(u_1-u_2)(q^{-2}u_1-u_2)(q^{-2}u_1-q^2u_1)v_L,
&if $L$ has shape $((1^2),\emptyset)$, \cr
0, &otherwise, \cr
}$$
and, when $u_1=0$,
$$pv_L = \cases{
(0-u_2)(-q^{-1}-q)(0-u_2)(0-u_2)v_L,
&if $L$ has shape $((1^2),\emptyset)$, \cr
0, &otherwise. \cr
}$$
Thus $p$ is an element of the ideal $I^{((1^2),\emptyset)}$.  Since
$I^{((1^2),\emptyset)}$ is a minimal ideal it is generated by
any one of its (nonzero) elements.
\endpf

\eject
\subsection The algebra $A_k(u_1,u_2;q)$

Let $u_1\in \CC$ and  $u_2, q\in \CC^*$.  Let
$A_k(u_1,u_2;q)$ be the algebra given by generators
$$X_1 \qquad\hbox{and}\qquad T_1,T_2,\ldots, T_{k-1}$$
and relations
$$
\matrix{
{\rm (1)}  & T_iT_j=T_jT_i,\hfill & |i-j|>1,\hfill \cr
{\rm (2)} & T_iT_{i+1}T_i=T_{i+1}T_iT_{i+1},\hfill & 1\le i\le n-2,
\hfill\cr
{\rm (3)} &T_i^2=(q-q^{-1})T_i+q, \hfill & 1\le i\le k-1,\hfill \cr
{\rm (4)} &X_1T_1X_1T_1=T_1X_1T_1X_1, \hfill\cr
{\rm (5)} &(X_1-u_1)(X_1-u_2)=0, \hfill \cr
{\rm (6)} &(X_1-u_2)(X_2-u_2)(X_2-q^2u_1)=0, \hfill
&\hbox{if $u_1\ne 0$,} \hfill \cr
\phantom{{\rm (6)}}
&(X_1-u_2)(T_1-q)(X_1-u_2)(X_2-u_2),\hfill
&\hbox{if $u_1=0$,\qquad where $X_2=T_1X_1T_1$.} \hfill \cr
\cr}\hskip.5truein
$$

Let
$$\hat A_k
= \{ (\lambda^{(1)},\lambda^{(2)}) \in \hat H_k^{(2)}
\ |\ \hbox{$\lambda^{(1)}$ has at most one row}\}.
\formula$$

\thm  Assume $H_k(u_1,u_2;q)$ is semisimple.
\smallskip\noindent
\item{(a)}  $A_k(u_1,u_2;q)$ is semisimple.
\smallskip\noindent
\item{(b)}
As in (2.14) let
$\langle I^{((1^2),\emptyset)}\rangle_k$ be the ideal of
$H_k(u_1,u_2;q)$
generated by the minimal ideal $I^{((1^2),\emptyset)}$ of
$H_2(u_1,u_2;q)$.  Then
$$A_k(u_1,u_2;q)
\cong {H_k(u_1,u_2;q)\over \langle I^{((1^2),\emptyset)}\rangle_k }. $$
\noindent
\item{(c)}
As in (2.13) let
$d_\lambda$ denote the number of standard tableaux of
shape $\lambda = (\lambda^{(1)},\lambda^{(2)})$ and
$M_d(\CC)$ the algebra of $d\times d$ matrices with
entries from $\CC$.  Then
$$A_k(u_1,u_2;q)\cong
\bigoplus_{\lambda\in \hat A_k} M_{d_\lambda}(\CC).$$
\smallskip\noindent
\item{(d)}  The irreducible $A_k(u_1,u_2;q)$-modules $H^\lambda$,
$\lambda\in
\hat A_k$,
are given by Hoefsmit's construction (Theorem 2.8).
\smallskip\noindent
\item{(e)} The  Bratelli diagram for the sequence of algebras
$A_1(u_1,u_2;q)\subseteq A_2(u_1,u_2;q)\subseteq \cdots \subseteq
A_k(u_1,u_2;q)$ has
vertices on level $m$ indexed by $\lambda\in \hat A_m$
and
edges $\lambda\longleftrightarrow \lambda^-$
if $\lambda^-$ is obtained from $\lambda$
by removing a box.  See Figure 2.
\pf
(a) follows from the fact that $A_k(u_1,u_2;q)$ is a quotient of
$H_k(u_1,u_2;q)$.
\smallskip\noindent
(b) By Lemma 2.15, the element $p$ generates the ideal
$\langle I^{((1^2),\emptyset)}\rangle$ in $H_k(u_1,u_2;q)$
and so this is a consequence of the definition of $A_k(u_1,u_2;q)$.
\smallskip\noindent
(c) By (2.14)
$$\langle I^{((1^2),\emptyset)}\rangle_k
= \bigoplus_{\lambda\supseteq ((1^2),\emptyset)} I^\lambda,
\formula$$
and so, by (b) and (2.13), the simple components of
$A_k(u_1,u_2;q)$ are indexed by those elements of
$\lambda\in \hat H_k^{(2)}$ which do not
contain $((1^2),\emptyset)$.
These are exactly the elements of $\hat A_k$.
\smallskip\noindent
(d) and (e) are conquences of (b), (c) and (2.18).
\endpf

$$ {\beginpicture
\setcoordinatesystem units <0.2cm,0.2cm>     
\setplotarea x from -3 to 33, y from 11 to 30     
\linethickness=0.05pt      

\put{$(\emptyset,\emptyset)$} at 20 30

\put{$(\emptyset, \phantom{t,,} )$} at 12 25
\put{$(\phantom{t,,} ,\emptyset)$} at 28 25
\putrectangle corners at 12.5 24.5 and 13.5 25.5
\putrectangle corners at 26.5 24.5 and 27.5 25.5
\plot 19 28.5 12 26.5 /
\plot 21 28.5 28 26.5 /

\put{$(\emptyset,\phantom{ttt,}  )$} at 4 20
\putrectangle corners at 4 19.5 and 5 20.5
\putrectangle corners at 6 19.5 and 5 20.5

\put{$(\emptyset,\phantom{tt,})$} at 12 20
\putrectangle corners at 12.5 21 and 13.5 20
\putrectangle corners at 12.5 19 and 13.5 20

\put{$(\phantom{tt},\phantom{tt})$} at 20 20
\putrectangle corners at 20.3 19.5 and 21.3 20.5
\putrectangle corners at 18.5 19.5 and 19.5 20.5

\put{$(\phantom{ttt,} ,  \emptyset)$} at 28 20
\putrectangle corners at 26.8 19.5 and 27.8 20.5
\putrectangle corners at 25.8 19.5 and 26.8 20.5

\plot  11 23.5  4 21.5 /
\plot  12 23.5  12 21.5 /
\plot  13 23.5  19 21.5 /
\plot  27 23.5  21 21.5 /
\plot  28 23.5  28 21.5 /

\put{$(\! \emptyset,\phantom{ttt,}\, )$} at -4.2 14
\putrectangle corners at -5 13.5 and -4 14.5
\putrectangle corners at -3 13.5 and -2 14.5
\putrectangle corners at -4 13.5 and -3 14.5

\put{$(\! \emptyset,\,\phantom{ttt})$} at 2.2 14
\putrectangle corners at 3 12.5 and 2 13.5
\putrectangle corners at 3 13.5 and 2 14.5
\putrectangle corners at 4 13.5 and 3 14.5

\put{$( \emptyset,\phantom{t}\, )$} at 7.3 14
\putrectangle corners at 8.5 13.5 and 7.5 14.5
\putrectangle corners at 8.5 13.5 and 7.5 12.5
\putrectangle corners at 8.5 11.5 and 7.5 12.5

\put{$( \phantom{t,}\,,\phantom{tt}\,)$} at 12.5 14
\putrectangle corners at 13.5 13.5 and 12.5 14.5
\putrectangle corners at 13.5 13.5 and 14.5 14.5
\putrectangle corners at 11.5 13.5 and 10.5 14.5

\put{$( \phantom{t,},\phantom{t,} )$} at 17.7 14
\putrectangle corners at 19 13.5 and 18 14.5
\putrectangle corners at 19 13.5 and 18 12.5
\putrectangle corners at 17.2 13.5 and 16.2 14.5

\put{$( \phantom{tt,}\, ,\phantom{t} )$} at 22.7 14
\putrectangle corners at 24.5 13.5 and 23.5 14.5
\putrectangle corners at 23 13.5 and 22 14.5
\putrectangle corners at 21 13.5 and 22 14.5

\put{$(\phantom{ttt,}\,\, ,  \emptyset\!)$} at 32.8 14
\putrectangle corners at 31.5 13.5 and 30.5 14.5
\putrectangle corners at 32.5 13.5 and 31.5 14.5
\putrectangle corners at 33.5 13.5 and 32.5 14.5

\plot 3 18.5  -3 15.5 /
\plot 4 18.5  2.5 15.5 /
\plot 5 18.5  12 15.5 /
\plot 19 18.5  13 15.5 /
\plot 20 18.5  18 15.5 /
\plot 11 18.5  3.5 15.5 /
\plot 11.5 18.5  8.5 15.5 /
\plot 13.5 18.5  17 15.5 /
\plot 20.5 18.5  22 15.5 /
\plot 27 18.5  23 15.5 /
\plot 28 18.5  32 15.5 /
\endpicture}
$$
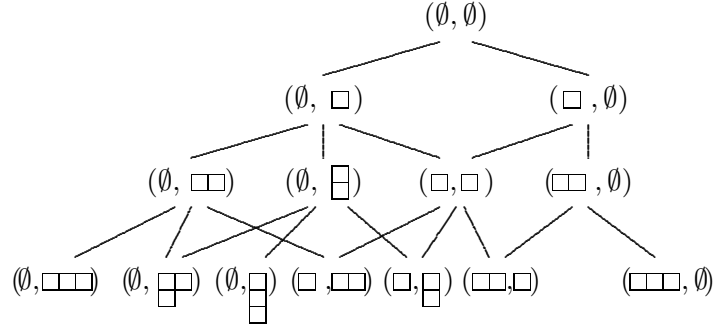
\centerline{Figure 2. Bratteli Diagram of $A_k(u_1,u_2;q)$ and $R_k(q)$}

\bigskip

By construction it is clear that the Bratteli diagram for
$A_k(u_1,u_2;q)$ is a subgraph of the Bratteli diagram for
$H_k(u_1,u_2;q)$.
It is the subgraph which is obtained by
removing all pairs of partitions
$\lambda=(\lambda^{(1)},\lambda^{(2)})$ which appear in (2.18)
(for all $k$).  Thus, it is the subgraph
which is obtained by removing the vertex $((1^2),\emptyset)$
and all its {\it descendants},
i.e., all pairs of partitions
$\lambda=(\lambda^{(1)}, \lambda^{(2)})$ which are
obtained by adding boxes to $((1^2),\emptyset)$.

\thm The algebra $A_k(u_1,u_2;q)$ is semisimple if and only if
$$\hbox{$q^{2d}u_1\ne u_2$ for all
$-k<d<k$,}\qquad\hbox{and}\qquad \quad [k]!\ne 0,$$
where $[k]! = [1][2]\cdots [k]$ and $[i] = 1+q^2+\cdots +q^{2(i-1)}$.
\pf
Since $A_k(u_1,u_2;q)$ is a quotient of $H_k(u_1,u_2;q)$, we know that
$A_k(u_1,u_2;q)$ is semisimple when $H_k(u_1,u_2;q)$ is.
Thus $A_k(u_1,u_2;q)$ is semisimple when (a) and (b) hold.

The Iwahori-Hecke algebra $H_k(q)$ of type $A_{k-1}$ is the cyclotomic
Hecke algebra $H_k(1,1;q)$. Thus,  $H_k(q)$ is the quotient
$A_k(u_1,u_2;q)$ by the relation $X_1=u_1$.  By Theorem 2.5
(in this case orginally due to Gyoja and Uno) $H_k(q)$ is semisimple
if and only if $[k]!\ne 0$.  Since $H_k(q)$ is a quotient
of $A_k(u_1,u_2;q)$, the algebra $A_k(u_1,u_2;q)$ is not semisimple
when $H_k(q)$ is not semisimple.  Thus, $A_k(u_1,u_2;q)$ is not
semisimple when $[k]!=0$.

If $u_1=u_2$ then the representation
$\rho\colon A_k(u_1,u_2;q)\to M_2(\CC)$ given by setting
$$\rho(X_1) = \pmatrix{ u_1 &1\cr 0 &u_1\cr}
\qquad\hbox{and}\qquad
\rho(T_i)=\pmatrix{q &0\cr 0 &q\cr}$$
is an indecomposable representation which is not irreducible.
Thus $A_k(u_1,u_1;q)$ is not semisimple.

The Tits deformation theorem (see [CR, (68.17)]) says that the
algebra $A_k(u_1,u_2;q)$ has the same structure for any choice
of the parameters $u_1,u_2$ for which it is semisimple.
Assume that $[k]!\ne 0$ and $u_2=q^{2d}u_1$, $u_1\ne 0$.
Let $\lambda/\mu$ be the skew shape given by
$\lambda = (k-1,d)$ and $\mu = (d-1)$ and define
$$CT(b) = u_1q^{2(c-r)+2},
\qquad\hbox{if $b$ is a box in position $(r,c)$ of $\lambda$.}$$
With these definitions the formulas in Theorem 2.2 define an
$H_k(u_1,q^{2d}u_1;q)$-module $H^{\lambda/\mu}$.  A check that
$$
 (X_1-q^{2d}u_1)(X_2-q^{2d}u_1)(X_2-q^2u_1)v_L=0,
$$
for all standard tableaux $L$ of shape $\lambda/\mu$ shows that
that $H^{\lambda/\mu}$ is an $A_k(u_1,q^{2d}u_1;q)$-module.
The standard proof (see [Ra, Theorem 4.1])
of Theorem 2.5 applies in this case to show that $H^{\lambda/\mu}$
is an irreducible $A_k(u_1,q^{2d}u_1;q)$-module.
It has dimension
$$\dim(H^{\lambda/\mu}) =\hbox{
(the number of standard tableaux of shape $\lambda/\mu$)}
= {k\choose d}-1.
$$
When we restrict this $A_k(u_1,q^{2d}u_1;q)$-module to $H_k(q)$,
it is a direct sum of irreducible $H_k(q)$-modules 
indexed by partitions $\nu\vdash k$ with multiplicity given 
by the classical Littlewood-Richardson coefficient 
$c_{\mu\nu}^\lambda$ (see [Ra, Theorem 6.1]). Since 
$\lambda = (k-1,d)$ and $\mu = (d-1)$ we have
$c_{\mu\nu}^\lambda\ne0$ only if $\nu$ has $\le 2$ rows
and $c_{\mu\nu}^\lambda = 0$ when $\nu = (k)$.
So $H^{\lambda/\mu}$ is an irreducible $A_k(u_1,q^2du_1;q)$-module
such that upon restriction to $H_k(q)$ is a direct sum of irreducible
representations indexed by partitions with length $\le 2$,
and which does not contain the ``trivial'' representation of $H_k(q)$.

When $A_k(u_1,u_2;q)$ is semisimple its irreducible representations 
$H^\lambda$ are indexed by pairs of parititions 
$\lambda=(\lambda^{(1)},\lambda^{(2)})$
such that $\lambda^{(1)}$ has at most one row.  If $\lambda^{(1)}$ has
length $r$, then, on restriction to $H_k(q)$, $H^\lambda$ is a direct 
sum of irreducibles indexed by partitions $\nu\vdash k$
with multiplicites $c_{\lambda^{(1)}\lambda^{(2)}}^\nu$.
By the Pieri rule [Mac, (5.16)] the resulting $\nu$ are 
those obtained by adding a horizontal 
strip of length $r$ to $\lambda^{(2)}$. Only when
$\lambda^{(2)}$ has a single row will all the $\nu$
have $\le 2$ rows  and, in this case, 
$c_{\lambda^{(1)}\lambda^{(2)}}^\nu=c_{(r),(k-r)}^\nu=1$ for $\nu=(k)$.
Thus, when $A_k(u_1,u_2;q)$ is semisimple
every irreducible representation which, on restriction
to $H_k(q)$, decomposes as a 
direct sum of components indexed by partitions with $\le 2$ 
rows does contain the ``trivial'' representation of $H_k(q)$.

Thus, the Tits deformation theorem implies that
$A_k(u_1,q^{2d}u_1;q)$ is not semisimple.
\endpf

\remark  It is interesting to note that the blob algebras 
(see [MW]) are also quotients of $A_k(u_1,u_2;q)$.

\subsection The $q$-rook monoid algebras $R_k(q)$

The new presentation of the $q$-rook monoid given in Section 1
shows that
$$
R_k(q) =A_k(0,1;q),
$$
and thus $R_k(q)$ is a quotient of Iwahori-Hecke algebra
$H_n(0,1;q)$ of type $B_k$.

\cor
\smallskip\noindent
(a)  The $q$-rook monoid algebra $R_k(q)$ is semisimple if and only if
$[k]_q!\ne 0$.
\smallskip\noindent
(b)  If $R_k(q)$ is semisimple then the irreducible representations
of $R_k(q)$ are indexed by $\lambda\in \hat A_k$ (see (2.16)) and are
given explicitly by the construction in Theorem 2.8.
\endthm

Part (a) of Corollary 2.21 is Theorem 2.19 applied
to $R_k(q)$ and (b) is Theorem 2.17(e) for $R_k(q)$.
Part (a) was proved in a different way by
Solomon [So4] and part (b) is the result
of Halverson [Ha, Theorem 3.2] which was the catalyst
for the results of this paper.

As in (2.10) it follows that, as $R_{k-1}(q)$-modules
$${\rm Res}^{R_k}_{R_{k-1}} R^\lambda = \bigoplus_{\lambda^-}
R^{\lambda^-},$$
where the sum runs over all pairs of partitions $\lambda^-$ which are
obtained
from $\lambda$ by removing a single box, and
$$R^{\lambda^-} = \CC\hbox{-span}\{ v_L\ |\ \hbox{$L^-$ has shape
$\lambda^-$}\},$$
where $L^-$ is the standard tableau with $k-1$ boxes which is obtained
by
removing the $k$ from $L$.
The first few rows of the Bratteli diagram for the sequence of algebras
$$R_1(q)\subseteq R_2(q)\subseteq R_3(q)\subseteq \cdots$$ are as
displayed in Figure 2.

\section 3.  Schur-Weyl dualities

Let $U_q\fgl(n)$ be the quantum group corresponding to
$GL_n(\CC)$.  This is the algebra given by generators
$$E_i,\quad F_i,\qquad (1\le i< n),
\qquad\hbox{and}\qquad q^{\pm\varepsilon_i},
\quad (1\le i\le n),$$
with relations
$$\matrix{
q^{\varepsilon_i}q^{\varepsilon_j} =
q^{\varepsilon_j}q^{\varepsilon_i},\quad
q^{\varepsilon_i}q^{-\varepsilon_i} =
q^{-\varepsilon_i}q^{\varepsilon_i}=1, \cr
\cr
q^{\varepsilon_i}e_jq^{-\varepsilon_i} =
\cases{ q^{-1}e_j, &if $j=i-1$, \cr
qe_j, &if $j=i$, \cr
e_j, &otherwise, \cr
}
\qquad\qquad
q^{\varepsilon_i}f_jq^{-\varepsilon_i} =
\cases{ qf_j, &if $j=i-1$, \cr
q^{-1}f_j, &if $j=i$, \cr
f_j, &otherwise, \cr
}\cr
\cr
\cr
\displaystyle{e_if_j-f_je_i = \delta_{ij}\,
{q^{\varepsilon_i-\varepsilon_{i+1}} -
q^{-(\varepsilon_i-\varepsilon_{i+1})}\over q-q^{-1}}, } \cr
\cr
\cr
e_{i\pm1}e_i^2-(q+q^{-1})e_ie_{i\pm1}e_i+e_i^2e_{i\pm1} =0,
\qquad
f_{i\pm1}f_i^2-(q+q^{-1})f_if_{i\pm1}f_i+f_i^2f_{i\pm1} =0, \cr
\cr
\cr
e_ie_j=e_je_i,\qquad\qquad f_if_j=f_jf_i, \qquad \hbox{if $|i-j|>1$}.
\cr
}$$
Part of the data of a quantum group is an ${\cal R}$-matrix, which
provides a canonical $U_q\fgl(n)$-module isomorphism
$$\check R_{MN}\colon M\otimes N \longrightarrow N\otimes M$$
for any two $U_q\fgl(n)$-modules $M$ and $N$.

The irreducible polynomial representations $L(\lambda)$ of
$U_q\fgl(n)$ are indexed by dominant integral weights
$\lambda\in L^+$ where
$$L = \sum_{i=1}^n \ZZ\varepsilon_i=\{\lambda=
\lambda_1\varepsilon_1+\cdots+\lambda_n\varepsilon_n
\ |\ \lambda_i\in \ZZ\},$$
and $L^+ = \{ \lambda\in L \ |\
\lambda_1\ge \lambda_2\ge \cdots\ge \lambda_n\}.$
The elements of $L^+$ can be identified with partitions
$\lambda$ with $\le n$ rows.

The irreducible representation
$$V=L(\varepsilon_1)
\qquad\hbox{has}\qquad
\dim(V) = n,
\qquad\hbox{and}\qquad
L(\mu)\otimes V \cong \bigoplus_{\mu^+} L(\mu^+)
\formula$$
as $U_q\fgl(n)$-modules, where the direct sum is over all partitions
$\mu^+$ which are obtained from $\mu$ by adding a box.
The $U_q\fgl(n)$-module $V$ can be given explicitly as the vector space
$$V = \hbox{$\CC$-span}\{v_1,\ldots, v_n\}$$
(so that the symbols $v_i$ form a basis of $V$) with
$U_q{\frak{gl}}(n)$-action given by
$$\matrix{
e_iv_j = \cases{ v_{j-1}, &if $j=i+1$, \cr 0, &if $j\ne i+1$, \cr}
\qquad
f_iv_j = \cases{ v_{j+1}, &if $j=i$, \cr 0, &if $j\ne i$, \cr}
\qquad\hbox{and} \cr
\cr
q^{\pm\varepsilon_i}v_j = \cases{ q^{\pm1}v_j, &if $j=i$, \cr v_j, &if
$j\ne i$.\cr}
\cr} $$
With this notation the $R$-matrix for $V\otimes V$
is given explicitly by
$$\eqalign{
&\check R_{VV}\colon V\otimes V\to V\otimes V,
\qquad\hbox{where}\qquad \cr
&\check R_{VV}(v_i\otimes v_j) = \cases{
qv_j\otimes v_i, &if $i=j$, \cr
v_j\otimes v_i, &if $i>j$, \cr
v_j\otimes v_i+(q-q^{-1})(v_i\otimes v_j), &if $i<j$. \cr
}}\formula$$

\subsection A Schur-Weyl duality for affine and cyclotomic Hecke
algebras

\thm (see [OR, Theorem 6.17ab and Theorem 6.18])
\smallskip\noindent
(a)  For any $\mu\in L^+$ there is an action of the
affine Hecke algebra
$\tilde H_k$ on $L(\mu)\otimes V^{\otimes k}$ given by
$\Phi\colon
\tilde H_k \longrightarrow \End(L(\mu)\otimes V^{\otimes k})$
where
$$
\Phi(X_1)
= \check R_{V,L(\mu)}\check R_{L(\mu),V}
\otimes \id_V^{\otimes (k-1)}
\qquad \hbox{and}\qquad
\Phi(T_i) = \id_{L(\mu)}\otimes \id_V^{\otimes(i-1)}\otimes
\check R_{VV} \otimes \id_V^{\otimes (k-i-1)}.
$$
(b)  The $\tilde H_k$ action on $L(\mu)\otimes V^{\otimes k}$
commutes with the $U_q\fgl(n)$-action and the map
$$\Phi\colon \tilde H_k \longrightarrow
\End_{U_q\fgl(n)}(L(\mu)\otimes V^{\otimes k})
\qquad\hbox{is surjective.}$$
\smallskip\noindent
(c) As a $(U_q\fgl(n),\tilde H_k)$ bimodule
$$L(\mu)\otimes V^{\otimes k} \cong
\bigoplus_\lambda L(\lambda)\otimes H^{\lambda/\mu},$$
where the sum is over all partitions $\lambda$ which are obtained from
$\mu$ by adding $k$ boxes and $H^{\lambda/\mu}$ is a simple
$\tilde
H_k$-module.
\smallskip\noindent
(d)  The representation $\Phi$ given in part (a) is a representation
of the cyclotomic Hecke algebra $H_k(u_1,\ldots, u_r;q)$, i.e.
$$\Phi\colon H_k(u_1,\ldots, u_r;q) \longrightarrow
\End_{U_q\fgl(n)}(L(\mu)\otimes V^{\otimes k}),$$
for any (multi)set of parameters $u_1,\ldots, u_r$
containing the (multi)set of values $CT(b)$ (defined in (2.1))
as $b$ runs over the boxes which can be added to $\mu$
(to get a partition).
\endthm

\remark  The affine Hecke algebra module $H^{\lambda/\mu}$ which
appears in Theorem 3.3(c) is the same as the module $H^{\lambda/\mu}$
constructed in Theorem 2.2.

\eject
\subsection A Schur-Weyl duality for Iwahori-Hecke algebras of type B

Suppose that $\mu$ is a partition with two addable boxes, i.e.
$$\mu =~~
{}_d\,\Bigg\lbrace
\underbrace{\beginpicture
\setcoordinatesystem units <0.28cm,0.28cm>     
\setplotarea x from 0 to 6, y from -1.5  to 2     
\linethickness=0.5pt     
\putrectangle corners at 0 2 and 6 -1.5
\endpicture}_{\ell} ~~= \ell^d,
\qquad\hbox{ for some $0<d\le n$, $\ell\in \ZZ_{\ge 0}$.}
$$

\prop  Let $\mu=\ell^d$, $0<d\le n$, $\ell\in \ZZ_{\ge 0}$.  Then
Theorem 3.3(d) provides a Schur-Weyl duality for $H_k(u_1,u_2;q)$
with $u_1=q^{2\ell}$ and $u_2=q^{2d}$.
\pf
One only needs to note that
if $b_1,b_2$ are the addable boxes
of $\mu$ and $CT(b)$ is as defined in (2.1) then
$$u_1 = CT(b_1)=q^{2\ell},
\qquad\hbox{and}\qquad u_2=CT(b_2)=q^{2d}. \qquad\qed$$

\subsection A Schur-Weyl duality for $A_n(u_1,u_2;q)$

Keeping the notation of Proposition 3.5, consider the special case
$\ell=n-1$ and $d\ge k$, so that
$$\mu =~~
{}_{n-1}\Bigg\lbrace
\underbrace{\beginpicture
\setcoordinatesystem units <0.28cm,0.28cm>     
\setplotarea x from 0 to 6, y from -1.5  to 2     
\linethickness=0.5pt     
\putrectangle corners at 0 2 and 6 -1.5
\endpicture}_{\ell} ~~= \ell^{(n-1)},
\qquad\hbox{ for some $\ell\in \ZZ_{\ge 0}$.}
$$
Then, as a $(U_q\fgl(n),H_k(u_1,u_2;q)$ bimodule
$$L(\mu)\otimes V^{\otimes k} \cong \bigoplus_\lambda
L(\lambda)\otimes H^\lambda,
\qquad\hbox{where}\qquad \lambda =
{\beginpicture
\setcoordinatesystem units <0.28cm,0.28cm>     
\setplotarea x from 0 to 6, y from -1.5  to 2     
\linethickness=0.5pt     
\putrectangle corners at 0 3 and 6 -1
\putrectangle corners at 0 -1 and  4 -2
\put{$\mu$} at 2.5 1
\put{${}^{{}^{\lambda^{(1)}}}$} at 2.5 -1.9
\put{${}^{\lambda^{(2)}}$} at 7.5 1.5
\plot 6  3  10  3 /
\plot 10 3  10  2 /
\plot 9  2  10  2 /
\plot 9  2  9   1 /
\plot 8  1  9   1 /
\plot 8  1  8  0 /
\plot 6 0  8  0 /
\endpicture}\formula$$
and $H^\lambda$ is a simple $H_k(u_1,u_2;q)$ module indexed by a
pair of partitions $\lambda = (\lambda^{(1)},\lambda^{(2)})$ with
$k$ boxes total and such that $\lambda^{(1)}$ has at most one row.
The following result shows that, in this case, the Schur-Weyl
duality in Theorem 3.3 becomes a Schur-Weyl duality for the algebra
$$A_k(u_1,u_2;q), \qquad\hbox{where}\qquad
u_1 = q^{2(n-1)} \quad\hbox{and}\quad u_2=q^{2\ell}.$$

\prop  Let $\mu = (n-1)^d$.  Then the $\tilde H_k$ action on
$L(\mu)\otimes V^{\otimes k}$ which is given by Theorem 3.3
factors through the algebra $A_k(u_1,u_2;q)$, where
$u_1=q^{2(n-1)}$ and $u_2=q^{2\ell}$.
\pf
By Theorem 3.3(d), the $\tilde H_k$ action on
$L(\mu)\otimes V^{\otimes k}$
factors through the algebra $H_k(u_1,u_2;q)$ where
$u_1=q^{2(n-1)}$ and $u_2=q^{2\ell}$.
It remains to check that
$(X_1-u_2)(X_2-u_2)X_2-q^2u_1)=0$ as operators on
$L(\mu)\otimes V^{\otimes k}$.
To do this it is sufficient to show that
$(X_1-u_2)(X_2-u_2)X_2-q^2u_1)=0$ as operators on
$H^\lambda$
for each $H^\lambda$ which appears in the decomposition (3.6).
The $H_k(u_1,u_2;q)$-module $H^\lambda$ has basis indexed by the
standard tableaux $L$ of shape $\lambda$ and
$$
(X_1-u_2)(X_2-u_2)(X_2-q^2u_1)v_L
=(CT(L(1))-u_2)(CT(L(2))-u_2)(CT(L(2))-q^2u_1)v_L.$$
For each of the possible positions of the first two boxes of
$L$ at least one of the factors in the last product is $0$.
Thus $(X_1-u_2)(X_2-u_2)X_2-q^2u_1)=0$ as operators on
$H^\lambda$.
\endpf

\eject
\subsection Another Schur-Weyl duality for cyclotomic Hecke algebras

Let $m_1,\ldots, m_r$ be positive integers such that
$m_1+\cdots m_r=n$.  Then
$\fg_P = \fgl(m_1)\oplus\cdots\oplus \fgl(m_r)$
is a Lie subalgebra of $\fgl(n)$,
and correspondingly
$$U_P = U_q\fgl(m_1)\otimes \cdots \otimes U_q\fgl(m_r)
\quad\hbox{is a subalgebra of}\quad
U_q\fgl(n).$$
There is a corresponding decomposition of
the fundamental representation
$V$ of $U_q\fgl(n)$ as a $U_P$-module:
$$V = V_1\oplus \cdots \oplus V_r,
\qquad\hbox{where}\qquad \dim(V_j)=m_j,$$
and $V_j$ is the fundamental representation for $U_q\fgl(m_j)$.
$$\hbox{If}\quad v\in V_j,\qquad\hbox{we write}\qquad
\deg(v)=j$$
and say that $v$ is {\it homogeneous of degree $j$}.

Let $u_1,\ldots, u_r\in \CC$ and define
$$d\colon V\to V
\qquad\hbox{by}\qquad d(v)=u_jv,\quad\hbox{if $\deg(v)=j$.}
\formula$$
Recall the action of $\check R_{VV}\colon V\otimes V\to V\otimes V$ as
given in (3.2), define
$$\check S_{VV}\colon V\otimes V\to V\otimes V
\qquad \hbox{by}\qquad
\check S_{VV}(v\otimes w) = \cases{
\check R_{VV}(v\otimes w), &if $\deg(v)=\deg(w)$,\cr
w\otimes v, &if $\deg(v)\ne \deg(w)$, \cr}$$
(for homogeneous $v,w\in V$), and define
$d_i, R_i, S_i\in \End(V^{\otimes k})$ by
$$\eqalign{
d_i &= \id_V^{\otimes(i-1)}\otimes d\otimes \id_V^{(k-i)}, \qquad
\qquad\quad 1\le i\le k, \cr
\check R_i &= \id_V^{\otimes(i-1)}\otimes
\check R_{VV}\otimes \id_V^{(k-i-1)},
\qquad 1\le i\le k-1, \cr
\check S_i &= \id_V^{\otimes(i-1)}\otimes
\check S_{VV}\otimes \id_V^{(k-i-1)},
\qquad
1\le i\le k-1. \cr
}\formula$$

\thm (Sakamoto-Shoji [SS])
\smallskip\noindent
(a)  There is an action of $H_k(u_1,u_2,\ldots, u_r;q)$
on $V^{\otimes k}$ given by
$\Phi_P\colon H_k(u_1,\ldots, u_r;q)\to
\End(V^{\otimes k})$ where
$$\Phi_P(T_i) = \check R_i
\qquad\hbox{and}\qquad
\Phi_P(X_1) = \check R_1^{-1}\cdots \check R_k^{-1}
\check S_k\cdots \check S_1 d_1.$$
\smallskip\noindent
(b) The action of $H_k(u_1,\ldots,u_r;q)$ commutes with the action
of $U_P$ on $V^{\otimes k}$,
\smallskip\noindent
\qquad
i.e., $\Phi_P\colon H_k(u_1,\ldots,u_k;q)\to\End_{U_P}(V^{\otimes k})$.
\smallskip\noindent
(c) As a $(U_P,H_k(u_1,\ldots, u_r;q))$-bimodule
$$V^{\otimes k} \cong
\bigoplus_{\lambda=(\lambda^{(1)},\ldots, \lambda^{(r)})}
L_P(\lambda)\otimes H^\lambda,$$
where the sum is over all $r$-tuples
$\lambda = (\lambda^{(1)},\ldots, \lambda^{(r)})$ of partitions such
that
$\ell(\lambda^{(j)})\le m_j$, $L_P(\lambda)$ is the simple $U_P$-module
given by
$$L_P(\lambda)=L^{(1)}(\lambda^{(1)})\otimes \cdots\otimes
L^{(r)}(\lambda^{(r)}),$$
where $L^{(j)}(\lambda^{(j)})$ is the simple $U_q(\fgl(m_j))$-module
correspnding
to the partition $\lambda^{(j)}$, and $H^\lambda$ is a
(not necessarily simple) $H_k(u_1,\ldots, u_r;q)$-module.
\endthm

\remark  The $H_k(u_1,\ldots, u_k;q)$-module $H^\lambda$
appearing in Theorem 3.10(c) is simple whenever
$H_k(u_1,\ldots, u_k;q)$ is semisimple.
In that case $H^\lambda$ coincides with the
$H_k(u_1,\ldots, u_r;q)$-module constructed in Theorem 2.8.

\eject
\subsection A Schur-Weyl duality for $R_k(q)$

Consider the case of Theorem 3.10 when
$$r=2,\qquad m_1=1,\qquad m_2=n,
\qquad u_1=0, \qquad u_2=1.$$
Then
$$V=V_1\oplus V_2
\qquad\hbox{where}\qquad
V_1=\CC v_0,
\quad\hbox{and}\quad
V_2 = \hbox{$\CC$-span}\{v_1,\ldots, v_n\}.$$
Let us analyze the action of $H_k(0,1;q)$ on $V^{\otimes k}$
as given by Sakamoto and Shoji.  Since $dv_0=u_1v_0=0$,
$$\Phi_P(X_1)(v_0\otimes v_{i_2}\otimes \cdots \otimes v_{i_k}) = 0,$$
and, for $\ell>0$,
$$\eqalign{
\Phi_P(X_1)(v_\ell\otimes v_{i_2}\cdots\otimes v_{i_k})
&= \check R_1^{-1}\cdots \check R_{k-1}^{-1}\check S_{k-1}\cdots \check
S_1d_1 (v_\ell\otimes v_{i_2}\otimes\cdots\otimes v_{i_k}) \cr
&= \check R_1^{-1}\cdots \check R_{k-1}^{-1}\check S_{k-1}\cdots \check
S_1 (v_\ell\otimes v_{i_2}\otimes\cdots\otimes v_{i_k}) \cr
&= \check R_1^{-1}\cdots \check R_{k-1}^{-1}\check R_{k-1}\cdots \check
R_1 (v_\ell\otimes v_{i_2}\otimes\cdots\otimes v_{i_k}) \cr
&= v_\ell\otimes v_{i_2}\otimes\cdots\otimes v_{i_k},\cr
}$$
and thus
$
\Phi_P(X_1) = d_1.
$
This calculation shows that the Sakamoto-Shoji action of $H_k(0,1;q)$
coincides exactly with action for the Schur-Weyl duality for
the $q$-rook monoid algebra $R_k(q)$ in the form given
by Halverson [Ha, Corollary 6.3].

\section 4. References

\bigskip

\item{[AK]} {\smallcaps S.\ Ariki and K.\ Koike},
{\it A Hecke algebra of $(\ZZ/r\ZZ)\wr S_n$ and construction of its
irreducible representations}, Adv.\ Math.\ {\bf 106} (1994), 216--243.

\medskip
\item{[Ar]} {\smallcaps S.\ Ariki},
{\it On the semisimplicity of the Hecke algbera of
$(\ZZ/r\ZZ)\wr S_n$}, J.\ Algebra {\bf 169} (1994),
216--225.

\medskip
\item{[Ch]} {\smallcaps I.\ Cherednik},
{\it A new interpretation of Gelfand-Tzetlin bases},
Duke Math.\ J.\ {\bf 54} (1987), 563-577.

\medskip
\item{[CR]} {\smallcaps C.\ Curtis and I.\ Reiner}, Methods of
Representation Theory: With Applications to Finite Groups and
Orders, Volume II, Wiley, New York, 1987.

\medskip
\item{[GU]} {\smallcaps A.\ Gyoja and K.\ Uno},
{\it On the semisimplicity of Hecke algebras},
J.\ Math.\ Soc.\ Japan {\bf 41} (1989), 75--79.

\medskip
\item{[Ha]} {\smallcaps T.\ Halverson},
{\it Representations of the $q$-rook monoid},
preprint (2001).

\medskip
\item{[Ho]} {\smallcaps P.\ N.\ Hoefsmit},
{\it Representations of Hecke algebras of finite
groups with $BN$-pairs of classical type}, Thesis, Univ.\ of British
Columbia, 1974.

\medskip
\item{[Mac]} {\smallcaps I.\ G.\ Macdonald},
{\rm Symmetric Functions and Hall polynomials,}
Second edition, Oxford University Press, New York, 1995.

\medskip
\item{[MW]} {\smallcaps P.\ P.\ Martin and D.\ Woodcock}, {\it On the
structure of the blob algebra}, J. Algebra {\bf 225} (2000),
957--988.
\medskip
\item{[Mu1]}  {\smallcaps W.\ D.\ Munn},
{\it Matrix representations of semigroups},
Proc.\ Camb.\ Phil.\ Soc., {\bf 53} (1957), 5-12.

\medskip
\item{[Mu2]} {\smallcaps W.\ D.\ Munn},
{\it The characters of the symmetric inverse
semigroup},
Proc.\ Camb.\ Phil.\ Soc., {\bf 53} (1957), 13-18.

\medskip
\item{[OR]} {\smallcaps R.\ Orellana and A.\ Ram},
{\it Affine braids, Markov traces
and the category ${\cal O}$}, preprint (2001).

\medskip
\item{[Ra]} {\smallcaps A.\ Ram}, {\it Skew shape representations
are irreducible}, preprint (1998).

\medskip
\item{[So1]} {\smallcaps L.\ Solomon},
{\it The Bruhat decomposition, Tits system
and Iwahori ring for the monoid of matrices over a finite field,\/}
Geom.\ Dedicata {\bf 36} (1990), 15--49.

\medskip
\item{[So2]} {\smallcaps L.\ Solomon},
{\it Representations of the rook monoid\/},
J.\ Algebra, to appear.

\medskip
\item{[So3]} {\smallcaps L.\ Solomon}, Abstract No.\
900-16-169, Abstracts Presented to the American Math.\
Soc., Vol.\ 16, No.\ 2, Spring 1995.

\medskip
\item{[So4]} {\smallcaps L.\ Solomon},
{\it The Iwahori algebra of ${\bf M}_n({\bf F}_q)$,
a presentation and a representation on
tensor space}, preprint (2001).

\medskip
\item{[SS]} {\smallcaps S.\ Sakamoto and T.\ Shoji},
{\it Schur-Weyl reciprocity for Ariki-Koike algebras},
J.\ Algebra {\bf 221} (1999), 293-314.

\medskip
\item{[Yg]} {\smallcaps A.\ Young},
{\it On quantitative substitutional analysis VI},
Proc.\ London Math.\ Soc.\  {\bf 31} (1931), 253-289.

\end